\documentclass[a4paper]{amsart}
\usepackage{amssymb}
\usepackage{graphicx}
\usepackage{color}
\usepackage{upref}
\usepackage{url}
\usepackage{stmaryrd}

\usepackage{calc}
\usepackage{ifthen}
\newcounter{hours}\newcounter{minutes}
\newcommand\printtime{%
 \setcounter{hours}{\time/60}%
 \setcounter{minutes}{\time-\value{hours}*60}%
\ifthenelse{\value{hours}<10}{0\thehours}{\thehours}
\ifthenelse{\value{minutes}<10}{:0\theminutes}{:\theminutes}}

\renewcommand{\phi}{\varphi}
\renewcommand{\epsilon}{\varepsilon}
\renewcommand{\theta}{\vartheta}

\newcommand{\dott}{\, \cdot\,}

\newcommand{\BV}{\mathbf{BV}}
\renewcommand{\d}[1]{\mathinner{\mathrm{d}{#1}}}
\newcommand{\C}[1]{\mathbf{C^{#1}}}
\newcommand{\Cc}[1]{\mathbf{C_c^{#1}}}
\renewcommand{\L}[1]{\mathbf{L^#1}}
\newcommand{\reali}{{\mathbb{R}}}

\newcommand{\naturali}{{\mathbb{N}}}
\newcommand{\modulo}[1]{{\left|#1\right|}}
\newcommand{\norma}[1]{{\left\|#1\right\|}}

\newcommand{\tv}{\mathop\mathrm{TV}}
\newcommand{\sgn}{\mathop\mathrm{sgn}}

\newcommand{\jmp}[1]{\left\llbracket #1\right\rrbracket}

\newtheorem{theorem}{Theorem}[section]
\newtheorem{proposition}[theorem]{Proposition}
\newtheorem{lemma}[theorem]{Lemma}
\newtheorem{definition}[theorem]{Definition}

\numberwithin{equation}{section}
\allowdisplaybreaks
\begin{document}

\title[Isentropic Fluid Dynamics in a Curved Pipe]{Isentropic Fluid
  Dynamics in a Curved Pipe}

\author[Colombo]{Rinaldo M.~Colombo} \address[Colombo]{\newline INDAM
  Unit, University of Brescia, Via Branze 38, I--25123 Brescia, Italy}
\email[]{Rinaldo.Colombo@Ing.UniBs.It}
\urladdr{http://dm.ing.unibs.it/rinaldo/}

\author[Holden]{Helge Holden} \address[Holden]{\newline Department of
  Mathematical Sciences, Norwegian University of Science and
  Technology, NO--7491 Trondheim, Norway}
\email[]{holden@math.ntnu.no}
\urladdr{http://www.math.ntnu.no/{\textasciitilde}holden}

\date{\today}



\thanks{The present work was supported by the GNAMPA~2015 project
  \emph{Balance Laws in the Modeling of Physical, Biological and
    Industrial Processes}, by the fund for international cooperation
  of the University of Brescia, as well as by the Research Council of
  Norway.}


\begin{abstract}
  In this paper we study isentropic flow in a curved pipe. We focus on
  the consequences of the geometry of the pipe on the dynamics of the flow. More
  precisely, we present the solution of the general Cauchy problem for
  isentropic  fluid flow in an arbitrarily curved,
  piecewise smooth pipe.  We consider initial data in the subsonic
  regime, with small total variation about a stationary solution. The
  proof relies on the front-tracking method and is based
  on~\cite{AmadoriGosseGuerra}.
\end{abstract}

\maketitle

\section{Introduction}
\label{sec:intro}

Consider a pipe filled with a compressible fluid. The pipe section is
far smaller than its length. The pipe is not assumed to be
rectilinear. We propose below a modification to the usual isentropic Euler
equations that takes into account the pipe's geometry.

\begin{figure}[!h]
  \centering 
\begin{picture}(0,0)%
\includegraphics{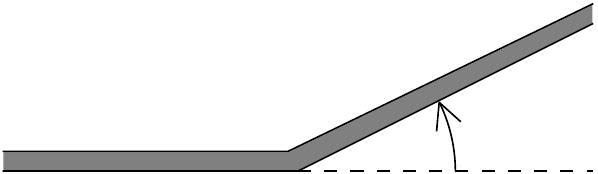}%
\end{picture}%
\setlength{\unitlength}{4144sp}%
\begingroup\makeatletter\ifx\SetFigFont\undefined%
\gdef\SetFigFont#1#2#3#4#5{%
  \reset@font\fontsize{#1}{#2pt}%
  \fontfamily{#3}\fontseries{#4}\fontshape{#5}%
  \selectfont}%
\fi\endgroup%
\begin{picture}(2724,789)(4489,-3673)
\put(6706,-3436){\makebox(0,0)[lb]{\smash{{\SetFigFont{12}{14.4}{\familydefault}{\mddefault}{\updefault}{\color[rgb]{0,0,0}$\vartheta$}%
}}}}
\end{picture}%
  \caption{A pipe curved by an angle $\theta$ at $x=0$, as considered
    in~\eqref{eq:3}--\eqref{eq:junction} or~\eqref{eq:deltaSource} and
    in Proposition~\ref{prop:RP}.}
  \label{fig:kink}
\end{figure}
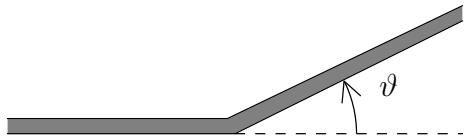

First, consider the case of a horizontal pipe with a single elbow.
Following~\cite{HoldenRisebroKink}, along the pipe we use the
classical isentropic $p$-system in Eulerian coordinates
\begin{equation}
  \label{eq:3}
  \left\{
    \begin{array}{l}
      \partial_t \rho + \partial_x q= 0,
      \\[2mm]
      \partial_t q + \partial_x \left(\dfrac{q^2}{\rho} + p (\rho)\right)
      =
      0,
    \end{array}
  \right.
\end{equation}
where $t$ is time, $x$ is the abscissa along the pipe, $\rho$ is the
mass density, $q$ is the linear momentum density, i.e., $q=\rho v$
where $v$ is the velocity, and $p$ is the pressure. At the kink,
located at, say, $x=0$, the following conditions on the traces of $q$
and of the dynamic pressure $P = q^2/\rho +p (\rho)$ are imposed:
\begin{equation}
  \label{eq:junction}
  q (t,0-) = q (t, 0+)
  \quad \mbox{ and } \quad
  P (t, 0-) = P (t, 0+) - f \,
  \kappa \left(2 \modulo{\sin (\theta/2)}\right) \, q (t, 0+) \,,
\end{equation}
where the positive parameter $f$ accounts for inhomogeneities in the
pipe's walls at the kink and $\kappa$ depends on the pipe's angle
$\theta$, see Figure~\ref{fig:kink}.
Equivalently, \eqref{eq:3}--\eqref{eq:junction} can be rephrased as a
single balance law with a Dirac delta source term in the second
equation:
\begin{equation}
  \label{eq:deltaSource}
  \left\{
    \begin{array}{l}
      \partial_t \rho + \partial_x q = 0,
      \\[2mm]
      \partial_t q + \partial_x \Big(\dfrac{q^2}{\rho} + p (\rho)\Big)
      =
      -f \, \kappa \left(2 \modulo{\sin (\theta/2)}\right) \, q \, \delta_{x=0}.
    \end{array}
  \right.
\end{equation}

Next, we consider a smoothly curved pipe described by the equation
$\Gamma = \Gamma(x)$ where $x$ is arc--length. It is reasonable to
assume that the dynamics of the fluid is governed by the equations
\begin{equation}
  \label{eq:2}
  \left\{
    \begin{array}{l}
      \partial_t \rho + \partial_x q = 0,
      \\[2mm]
      \partial_t q + \partial_x \left(\dfrac{q^2}{\rho} + p (\rho)\right)
      =
      -f (x) \; \kappa \left(\norma{\Gamma'' (x)}\right) \, q \,,
    \end{array}
  \right.
\end{equation}
where $\norma{\Gamma'' (x)}$ equals the curvature of the pipe at the
location $\Gamma(x)$.  We have $\kappa(0)=0$, and $f(x)$ is an empirical factor that depends
on the location along the pipe.  A brief derivation of the model can be found in \cite{HoldenRisebroKink}.

More generally, we consider an arbitrary piecewise smooth pipe. Call
$\bar x_0, \ldots, \bar x_m$ its corner points, or kinks, and denote by
$\theta_i$ the angle of the pipe at $\bar x_i$, see
Figure~\ref{fig:kink}.  To avoid unphysical behavior we assume that the pipe is horizontal and
rectilinear outside a compact set.  We
are thus led to consider the system:
\begin{equation}
  \label{eq:smoothG}
  \left\{
    \begin{array}{@{}l@{\;}l@{}}
      \partial_t \rho + \partial_x q = 0,
      \\[2mm]
      \partial_t q + \partial_x \left(\dfrac{q^2}{\rho} + p (\rho)\right)
     & =

      -f (x) \; \kappa \left(\norma{\Gamma'' (x)}\right) \, q
      - \rho \, g \, \sin \alpha (x)
      \\[1mm]
      &\quad
      \displaystyle
      -\sum_{i=0}^m
      f (\bar x_i) \, \kappa \left(2 \modulo{\sin (\theta_i/2)}\right)
      \, q (t, \bar x_i+) \, \delta_{x=\bar x_i} \,,
    \end{array}
  \right.
\end{equation}
where $\alpha = \alpha (x)$ describes the inclination of the pipe with respect to
the horizontal plane  at $x$ and $g$ is gravity.  Both $\kappa$ and $\alpha$ vanish outside a compact set.

The main result of the present paper is that~\eqref{eq:smoothG}
generates a Lipschitz continuous semigroup defined globally in time on
all initial data that are small perturbations of stationary
solutions. The results in~\cite{AmadoriGosseGuerra} also ensure the
uniqueness of this semigroup.

The analytic techniques employed here are rooted in the idea of
approximating the piecewise smooth pipe with a polygonal one. Indeed,
the case of a polygonal pipe can be obtained by \emph{gluing together}
systems of the type~\eqref{eq:junction}, where the source is a
sequence of linear combinations of Dirac delta masses, which
correspond to stationary discontinuities. At this point, the front-tracking
method for systems of hyperbolic conservation laws~\cite{BressanLectureNotes, HoldenRisebro}
proves to be a very effective tool.
First, front-tracking approximations are defined through the available solutions to
Riemann problems, including those at the Dirac masses. Second,
front-tracking approximations are extremely accurate in capturing the
essential features of the exact solutions to conservation laws. Third,
analytic techniques are available that allow to prove the convergence
of these approximations. We refer to~\cite{BressanLectureNotes,
  HoldenRisebro} for further details on the front-tracking method.

\section{Main Result}
\label{sec:m}

Throughout this paper, $\reali^+ = (0, +\infty)$ and
$\overline\reali^+ = [0, +\infty)$. Moreover, we denote the state of
the fluid by $u$, where $u \equiv (\rho,q)$, with $q = \rho \,v$.

We assume that the fluid can be described through the pressure law $p$
satisfying
\begin{description}
\item[($\boldsymbol{p}$)]$p \in \C2 (\reali^+; \reali^+)$, $p' (\rho)
  \geq 0$ and $p'' (\rho)\geq 0$ for all $\rho >0$.
\end{description}
\noindent A typical example is a polytropic gas with the $\gamma$-pressure law $p(\rho) =
\rho^\gamma$ for $\gamma \geq 1$.

With reference to the $p$-system~\eqref{eq:3} recall the following
quantities
\begin{equation}
  \label{eq:EFP}
  \begin{aligned}
    E (\rho,q) & = \frac{q^2}{2\rho} + \rho \int_{\bar\rho}^\rho
    \frac{p (r)}{r} \d{r}, & & \mbox{mathematical entropy,}
    \\
    F (\rho,q) & = \frac{q}{\rho} \left(E (\rho,q) + p (\rho)\right),
    & &\mbox{entropy flow,}
    \\
    P (\rho,q) & = \frac{q^2}{\rho} + p (\rho), & & \mbox{dynamic
      pressure.}
  \end{aligned}
\end{equation}

\subsection{Stationary Solutions}
\label{subsec:SS}

Assume the pipe is horizontal. Then, both
systems~\eqref{eq:deltaSource} and~\eqref{eq:2} admit the stationary
solution
\begin{displaymath}
  q = 0 \quad \mbox{ and } \quad \rho = \mbox{constant.}
\end{displaymath}

In the case of a single kink~\eqref{eq:deltaSource}, further
stationary solutions are given by
\begin{align*}
    \rho & =
      \begin{cases}
        \rho^\ell, & x  <  0,
        \\
        \rho^r,&   x  >  0,
      \end{cases}
\quad   q  =  \mbox{constant,  where }
  P (\rho^\ell, q) - P (\rho^r, q)
  = -f \kappa\left(2\modulo{\sin (\theta/2)}\right) q \,.
\end{align*}
Stationary solutions in the case of a polygonal pipe are obtained
by gluing together solutions of the type above, i.e., $q$ is constant
while $\rho$ satisfies the jump condition at every kink.

In a smooth pipe with gravity, stationary solutions satisfy
\begin{displaymath}
  \partial_x P\left(\rho (x), q\right)
  =
  -f (x) \, \kappa\left(\norma{\Gamma'' (x)}\right) \, q
  - \rho \, g \, \sin \alpha (x)
  \quad \mbox{ and } \quad
  q = \mbox{constant.}
\end{displaymath}

Gluing together stationary solutions of the types above yields stationary
solutions in the case of a piecewise smooth pipe.

Throughout this paper, by $\bar u = \bar u (x)$ we denote any of the
stationary solutions constructed above.

\subsection{The Case of a Single Kink}
\label{subsec:RP}

We now briefly consider the Riemann Problem
for~\eqref{eq:deltaSource}, referring to~\cite{HoldenRisebroKink} for
more details.

The pipe consists now of two rectilinear tubes connected through a
kink at an angle $\theta \in (- \pi, \pi)$ located at, say, $x=0$, so
that
\begin{displaymath}
  \Gamma (x) =
    \begin{cases}
      (1,0) \, x, & x  <  0,
      \\
      (\cos \theta, \sin \theta) \, x, & x  >  0,
    \end{cases}
\end{displaymath}
see Figure~\ref{fig:kink}.  Then, the Riemann Problem for the
model~\eqref{eq:3}--\eqref{eq:junction} or~\eqref{eq:deltaSource}
introduced in~\cite{HoldenRisebroKink} reads
\begin{equation}
  \label{eq:1}
  \left\{
    \begin{array}{l}
      \partial_t \rho + \partial_x q = 0,
      \\[2mm]
      \partial_t q
      +
      \partial_x P (\rho,q) = 0,
      \\[2mm]
      \jmp{q} (t,0)
      =
      0,
      \\[2mm]
      \jmp{P} (t,0)
      =
      f \, \kappa\left(2 \modulo{\sin (\theta/2)}\right) \, q (t, 0+),
      \\[2mm]
      (\rho,q) (0,x) =
        \begin{cases}
          (\rho^l,q^l), & x  <  0,
          \\[1mm]
          (\rho^r,q^r), & x  >  0,
        \end{cases}
        \end{array}
   \right.
\end{equation}
where, as usual, we denote\footnote{Here $F(x\pm)=\lim_{h\downarrow
    0}F(x\pm h)$ for any function $F$.}  $\jmp{F} (t,x) = F (t, x+) -
F (t, x-)$ for any function $F$ of the pair $(\rho,q)$. The function
$\kappa$ is assumed to satisfy
\begin{description}
\item[($\boldsymbol{\kappa}$)] $\kappa \in \C1 (\reali;
  \overline\reali^+)$, with $\kappa (0) = 0$ and $\kappa$ is even.
\end{description}
\noindent We also introduce the \emph{subsonic region}
\begin{equation}
  \label{eq:Omega}
  \Omega =
  \left\{
    (\rho,q) \in \reali^+ \times \reali \colon
    \modulo{q / \rho} < \sqrt{p' (\rho)}
  \right\} \,,
\end{equation}
where the velocity $v=q/\rho$ of the fluid is smaller than the sound speed $\sqrt{p' (\rho)}$.
 Due to its relevance in the
applications, we restrict our attention below to initial data and solutions
attaining values in the subsonic region.

\begin{proposition}
  \label{prop:RP}
  Let~\textbf{($\boldsymbol{p}$)} and~\textbf{($\boldsymbol{\kappa}$)}
  hold. Fix $f>0$ and a subsonic stationary solution $\bar u$
  to~\eqref{eq:1}. Then, there exists a $\delta > 0$ such that for all
  states $u^\ell, u^r \in \Omega$ satisfying
  \begin{displaymath}
    \norma{u_o - \bar u}_{\L\infty (\reali; \reali^+ \times \reali)} < \delta
    \quad \mbox{ where } \quad
    u_o (x) =
      \begin{cases}
        u^\ell, & x <  0 \,,
        \\u^r, & x  >  0 \,,
      \end{cases}
  \end{displaymath}
  the Riemann Problem~\eqref{eq:1} admits a unique self-similar weak
  entropy solution attaining values in $\Omega$, consisting of a
  $1$-wave supported in $x<0$, a jump along $x=0$ and a $2$-wave
  supported in $x>0$.
\end{proposition}

The Riemann Problem~\eqref{eq:1} was analyzed for arbitrary initial states
in~\cite[Section~2]{HoldenRisebroKink} in the isothermal case where the pressure $p
(\rho) = \rho$. The well known properties of the $p$-system allow us
to apply~\cite[Theorem~3.2]{ColomboHertySachers}, so that the Cauchy
problem for~\eqref{eq:1} is well posed in $\L1$. The proof of
Proposition~\ref{prop:RP} directly follows from the cited references.

\subsection{The Case of a Piecewise Smooth Pipe}
\label{subsec:smooth}

We now consider a piecewise smooth pipe with finite curvature, see Figure~\ref{fig:pipe}.  More precisely, we make the following assumptions:
\begin{figure}[!h]
  \centering 
\begin{picture}(0,0)%
\includegraphics{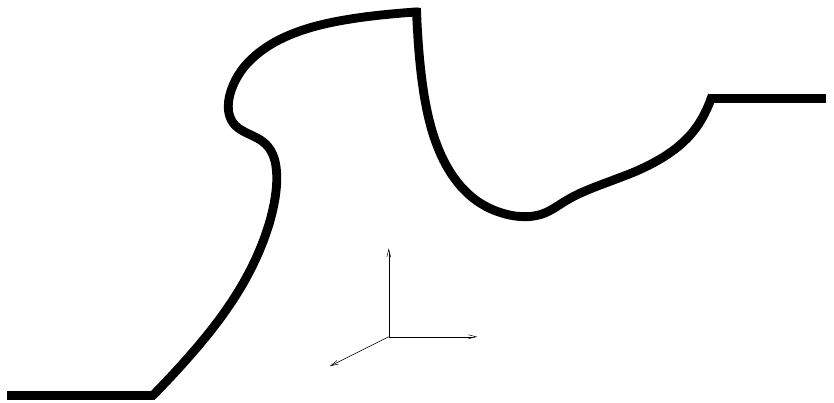}%
\end{picture}%
\setlength{\unitlength}{1243sp}%
\begingroup\makeatletter\ifx\SetFigFont\undefined%
\gdef\SetFigFont#1#2#3#4#5{%
  \reset@font\fontsize{#1}{#2pt}%
  \fontfamily{#3}\fontseries{#4}\fontshape{#5}%
  \selectfont}%
\fi\endgroup%
\begin{picture}(12694,6124)(-523,-7818)
\put(5626,-5911){\makebox(0,0)[lb]{\smash{{\SetFigFont{10}{6.0}{\familydefault}{\mddefault}{\updefault}{\color[rgb]{0,0,0}$\mathbf{k}$}%
}}}}
\put(6241,-2071){\makebox(0,0)[lb]{\smash{{\SetFigFont{10}{6.0}{\familydefault}{\mddefault}{\updefault}{\color[rgb]{0,0,0}$\Gamma(\bar x_i)$}%
}}}}
\end{picture}%
  \caption{A piecewise smooth pipe.}
  \label{fig:pipe}
\end{figure}
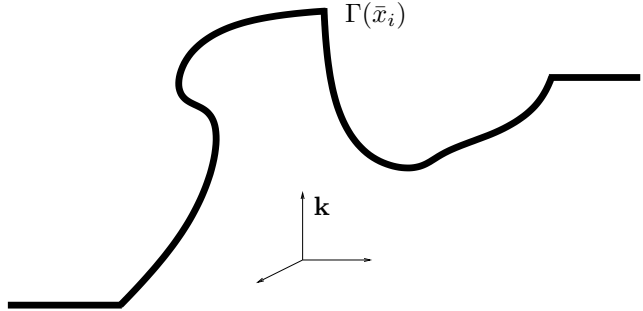
\begin{description}
\item[($\boldsymbol\Gamma$)] $\Gamma \in \C0 (\reali; \reali^3)$ is
  such that:
  \begin{enumerate}
  \item $\Gamma$ is piecewise smooth: there exist $\bar x_0, \bar x_1,
    \ldots, \bar x_m$ with $x_{i-1} < x_i$ for all $i$ such that
    $\Gamma|_{(-\infty, \bar x_0]} \in \C2 ((-\infty, \bar x_0];
    \reali^3)$, $\Gamma|_{[\bar x_{i-1}, \bar x_i]} \in \C2 ([\bar
    x_{i-1}, \bar x_i], \reali^3)$ and $\Gamma|_{[\bar x_m, +\infty)}
    \in \C2 ([\bar x_m, +\infty); \reali^3)$;
  \item $\Gamma$ is parametrized by arc--length: $\norma{\Gamma' (x)}
    = 1$ for all $x \in \reali\setminus\{\bar x_0, \ldots, \bar
    x_m\}$;
  \item $\Gamma$ has finite curvature: $\Gamma''$ vanishes outside a
    compact set;
  \item $\Gamma$ is horizontal outside a compact set: $\Gamma' (x)
    \cdot {\mathbf k}$ vanishes outside a compact, where ${\mathbf k}$ denotes the unit
    vertical vector.
  \end{enumerate}
\end{description}

On the friction term $f$, we require the following condition:
\begin{description}
\item[($\boldsymbol{f}$)] $f \in (\C0 \cap \L\infty) (\reali;\reali)$
  and $f\geq 0$.
\end{description}

\begin{lemma}
  Let $\Gamma$ satisfy~\textbf{($\mathbf{\Gamma}$)}. Then, $\Gamma'
  \in \BV (\reali; \reali^3)$ and its weak derivative is the measure
  \begin{displaymath}
    \mu
    =
    \Gamma'' \, \d{\mathcal{L}}
    +
    \sum_{i=0}^m \left(\Gamma' (x_i+) - \Gamma' (x_i-)\right) \delta_{x=x_i} \,.
  \end{displaymath}
\end{lemma}

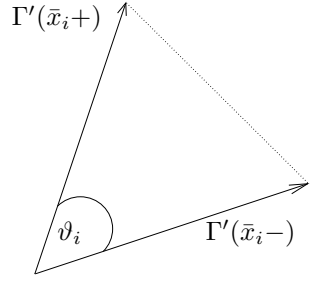
\begin{figure}[!h]
  \centering 
\begin{picture}(0,0)%
\includegraphics{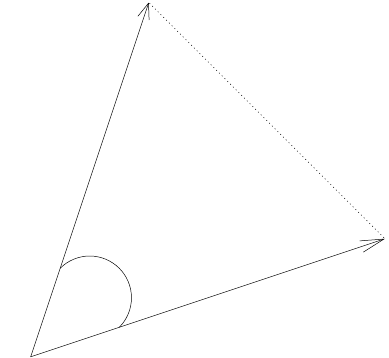}%
\end{picture}%
\setlength{\unitlength}{1243sp}%
\begingroup\makeatletter\ifx\SetFigFont\undefined%
\gdef\SetFigFont#1#2#3#4#5{%
  \reset@font\fontsize{#1}{#2pt}%
  \fontfamily{#3}\fontseries{#4}\fontshape{#5}%
  \selectfont}%
\fi\endgroup%
\begin{picture}(5877,5424)(1786,-5923)
\put(2701,-5236){\makebox(0,0)[lb]{\smash{{\SetFigFont{10}{6.0}{\familydefault}{\mddefault}{\updefault}{\color[rgb]{0,0,0}$\theta_i$}%
}}}}
\put(1801,-961){\makebox(0,0)[lb]{\smash{{\SetFigFont{10}{6.0}{\familydefault}{\mddefault}{\updefault}{\color[rgb]{0,0,0}$\Gamma'(\bar x_i+)$}%
}}}}
\put(5626,-5236){\makebox(0,0)[lb]{\smash{{\SetFigFont{10}{6.0}{\familydefault}{\mddefault}{\updefault}{\color[rgb]{0,0,0}$\Gamma'(\bar x_i-)$}%
}}}}
\end{picture}%
  \caption{Justification of~\eqref{eq:CompGamma}: here,
    $\norma{\Gamma' (\bar x_i+)} = \norma{\Gamma' (\bar x_i+)} = 1$.}
  \label{fig:Gamma}
\end{figure}
Remark that the above expression of $\mu$ admits a geometric
interpretation. For $i= 0, \ldots, m$, call $\theta_i$ the angle at
$\bar x_i$ such that $\cos \theta_i = \Gamma' (\bar x_i-) \cdot
\Gamma' (\bar x_i+)$. Elementary geometric considerations, see
Figure~\ref{fig:Gamma}, show that
\begin{equation}
  \label{eq:CompGamma}
  \norma{\Gamma' (\bar x_i+) - \Gamma' (\bar x_i-)}
  =
  \sqrt{2 (1-\cos \theta_i)}
  = 2 \modulo{\sin (\theta_i/2)},
\end{equation}
as used in~\cite{HoldenRisebroKink}.

\begin{definition}
  \label{def:solPlimit}
  Let $T>0$ and fix a stationary state $\bar u \in \reali^+ \times
  \reali$. By a \emph{weak solution} to~\eqref{eq:2} we mean a map
  \begin{displaymath}
    u = (\rho,q) \in
    \C0\left([0,T];
      \bar u
      +
      (\L1 \cap \BV) (\reali; \reali^+ \times \reali)
    \right)
  \end{displaymath}
  such that $u_o = u_{|t=0}$ and for any function $\phi \in \Cc1 ((0,
  T) \times \reali; \reali)$, we have
  \begin{align*}
    \int_\reali \int_0^T \left( \rho \, \partial_t \phi + q
      \, \partial_x \phi \right) \d{t} \, \d{x} =& 0,
    \\
    \int_\reali \int_0^T \left( q \, \partial_t \phi + P
      (\rho,q) \partial_x \phi \right) \d{t} \, \d{x} = &
    \int_{\reali} \int_0^T f (x) \, \kappa \! \left(\norma{\Gamma''
        (x)}\right) \, q (t,x) \, \phi (t,x) \d{t} \, \d{x}
    \\
    & + \sum_{i=0}^m \int_0^T f (\bar x_i) \, \kappa\left(2 \sin
      (\theta_i/2)\right) \, q (t, \bar x_i) \, \phi(t, \bar x_i)
    \d{t}
    \\
    & + \int_{\reali} \int_0^T \rho (t,x) \, g \, \sin \alpha (x) \,
    \phi (t,x) \d{t} \d{x}.
  \end{align*}
  The weak solution $(\rho,q)$ is a \emph{weak entropy solution} if
  for any function $\phi \in \Cc1 ((0, T) \times \reali; \reali^+)$,
  we have
  \begin{align*}
 &   \int_{\reali} \int_0^T
    \left(
      E (\rho,q) \, \partial_t \phi
      +
      F (\rho,q) \, \partial_x \phi
    \right) \d{t} \d{x}
   \\
&    +
    \int_{\reali} \int_0^T
      \partial_q E (\rho,q) \big(
        f (x) \, \kappa \! \left(\norma{\Gamma''
            (x)}\right) \, q (t,x)
        +
        \rho (t,x) \, g \, \sin \alpha (x)
      \big) \phi
 \d{t} \d{x}
     \geq  0 \,.
  \end{align*}
\end{definition}

\begin{theorem}
  \label{thm:limit}
  Let~\textbf{($\boldsymbol{p}$)}, \textbf{($\boldsymbol{\Gamma}$)},
  \textbf{($\boldsymbol{f}$)}, and~\textbf{($\boldsymbol{\kappa}$)}
  hold. Fix a subsonic stationary solution $\bar u$.  Then, there
  exist $\hat\delta, \check\delta$, and $L \in \reali^+$ such
  that~\eqref{eq:2} generates a semigroup
  \begin{displaymath}
    S \colon \reali^+ \times \mathcal{D} \to \mathcal{D}
  \end{displaymath}
  with the properties:
  \begin{enumerate}
  \item The domain $\mathcal{D}$ is non-trivial and its elements have
    uniformly bounded total variation:
    \begin{align*}
      \left\{ u \in \bar u + \L1 (\reali;\reali^+ \times \reali)
        \colon \tv (u) \leq \check \delta \right\} & \subseteq
      \mathcal{D},
      \\
      \left\{ u \in \bar u + \L1 (\reali;\reali^+ \times \reali)
        \colon \tv (u) \leq \hat \delta \right\} & \supseteq
      \mathcal{D}.
    \end{align*}
  \item For all $u_o \in \mathcal{D}$, the map $t \mapsto S_t u_o$ is
    a weak entropy solution to~\eqref{eq:smoothG} in the sense of
    Definition~\ref{def:solPlimit}.
  \item $S$ is Lipschitz continuous with respect to the $\L1$ norm,
    i.e., for $u, u' \in \mathcal{D}$
    \begin{displaymath}
      \norma{S_{t'} u' - S_{t} u}_{\L1 (\reali;\reali^+ \times \reali)}
      \leq L \left(
        \norma{u' - u}_{\L1 (\reali;\reali^+ \times \reali)}
        +
        \modulo{t'-t}
      \right).
    \end{displaymath}
  \end{enumerate}
\end{theorem}

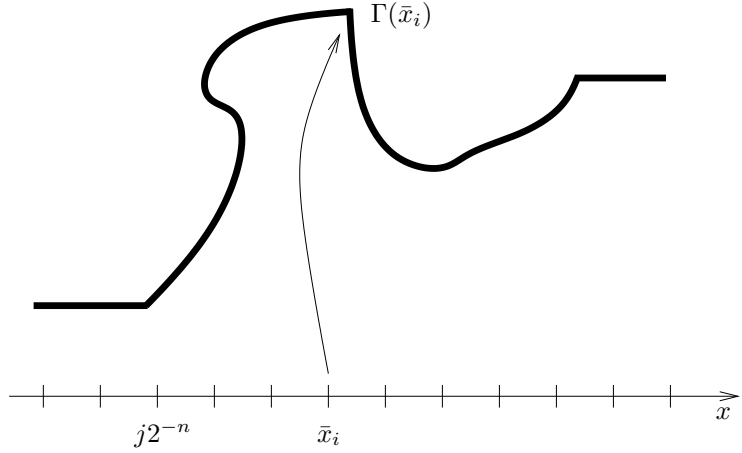
\begin{figure}[!h]
  \centering 
\begin{picture}(0,0)%
\includegraphics{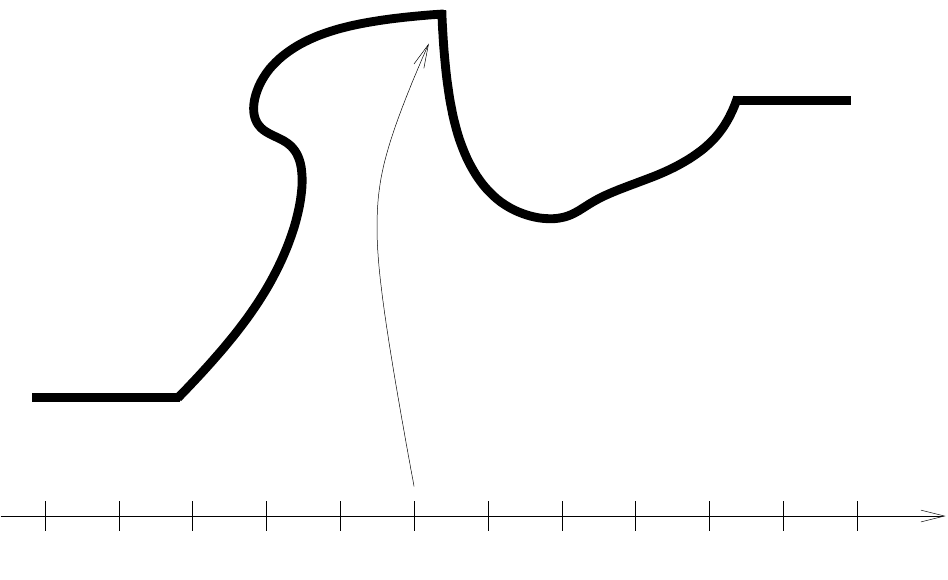}%
\end{picture}%
\setlength{\unitlength}{1243sp}%
\begingroup\makeatletter\ifx\SetFigFont\undefined%
\gdef\SetFigFont#1#2#3#4#5{%
  \reset@font\fontsize{#1}{#2pt}%
  \fontfamily{#3}\fontseries{#4}\fontshape{#5}%
  \selectfont}%
\fi\endgroup%
\begin{picture}(14424,8790)(-911,-10484)
\put(6241,-2071){\makebox(0,0)[lb]{\smash{{\SetFigFont{10}{6.0}{\familydefault}{\mddefault}{\updefault}{\color[rgb]{0,0,0}$\Gamma(\bar x_i)$}%
}}}}
\put(13051,-9961){\makebox(0,0)[lb]{\smash{{\SetFigFont{10}{6.0}{\familydefault}{\mddefault}{\updefault}{\color[rgb]{0,0,0}$x$}%
}}}}
\put(5176,-10411){\makebox(0,0)[lb]{\smash{{\SetFigFont{10}{6.0}{\familydefault}{\mddefault}{\updefault}{\color[rgb]{0,0,0}$\bar x_i$}%
}}}}
\put(1576,-10411){\makebox(0,0)[lb]{\smash{{\SetFigFont{10}{6.0}{\familydefault}{\mddefault}{\updefault}{\color[rgb]{0,0,0}$j 2^{-n}$}%
}}}}
\end{picture}%
  \caption{Discretization of~\eqref{eq:smoothG} leading
    to~\eqref{eq:delta}.}
  \label{fig:pipe2}
\end{figure}

\begin{proof}
  We follow the construction in~\cite{AmadoriGosseGuerra}. In the
  discretization of the pipe, we assume for simplicity that all kinks
  are at a dyadic abscissa. In other words, without any loss of
  generality, we assume that for all $i = 0, \ldots, m$, we have $\bar
  x_i = j_i 2^{-n_i}$ for suitable $n_i \in \naturali$ and $j_i \in
  \{-2^{2n_i}, \ldots, 2^{2n_i}\}$, see Figure~\ref{fig:pipe2}.

  Introduce the set $\mathcal{K}_n$ of indices that correspond to
  kinks, namely
  \begin{displaymath}
    \mathcal{K}_n
    =
    \left\{
      j \in \{-2^{2n}, \ldots, 2^{2n}\} \colon
      \exists i \in \{0, \ldots, m\} \mbox{ such that }
      \bar x_i = j \, 2^{-n}
    \right\} \,.
  \end{displaymath}
  The procedure in~\cite[Theorem~3]{AmadoriGosseGuerra}, by means of
  front-tracking approximate solutions to~\eqref{eq:smoothG},
  constructs an \emph{exact} solution $u^n$ to the following
  approximation of~\eqref{eq:smoothG}:
  \begin{equation}
    \label{eq:delta}
    \left\{
      \begin{array}{@{}ll@{}}
        \partial_t \rho + \partial_x q = 0,
        \\[2mm]
        \partial_t q + \partial_x \left(\dfrac{q^2}{\rho} + p (\rho)\right)
        =
        &\!\!\!\!
        \displaystyle
        - \!\!
        \sum_{j \not\in \mathcal{K}_n} \!\!
        f (j2^{-n}) \, \kappa \left(\norma{\Gamma'' (j2^{-n})}\right)
        q(t,j2^{-n}) \, \delta_{x=j2^{-n}}
        \\[4mm]
        &\!\!\!\!
        \displaystyle
        - \!\!
        \sum_{j =-2^{2n}}^{2^{2n}}
        \rho(t,j2^{-n}) \, g \, \sin \alpha (j2^{-n}) \, \delta_{x=j2^{-n}}
        \\[4mm]
        &\!\!\!\!
        \displaystyle
        - \!\!
        \sum_{i=0}^m
        f (\bar x_i) \, \kappa \left(2 \modulo{\sin (\theta_i/2)}\right)
        \, q (t, \bar x_i+) \, \delta_{x=\bar x_i} \,.
      \end{array}
    \right.\!\!\!
  \end{equation}
  An application of~\cite[Theorem~6]{AmadoriGosseGuerra} yields for
  any $n \in \naturali$ the existence of a semigroup $S^n \colon
  \reali^+ \times \mathcal{D}^n \to \mathcal{D}^n$ satisfying~(1) with
  $\mathcal{D}$ replaced by $\mathcal{D}^n$, (2)
  with~\eqref{eq:smoothG} replaced by~\eqref{eq:delta}, and~(3) for
  suitable $\hat\delta, \check\delta$ and $L$ independent of $n$.

  We now let $n \to +\infty$ and follow the procedure
  in~\cite[Theorem~8]{AmadoriGosseGuerra}. Remark
  that~\eqref{eq:delta} differs from the equation considered
  in~\cite{AmadoriGosseGuerra} by the last term
  \begin{displaymath}
    -
    \sum_{i=0}^m
    f (\bar x_i) \, \kappa \left(2 \modulo{\sin (\theta_i/2)}\right)
    \, q (t, \bar x_i+) \, \delta_{x=\bar x_i}
  \end{displaymath}
  on the right-hand side of the second equation. However, this term is
  independent of $n$ and does not prevent the application of
  techniques used in~\cite[Theorem~8]{AmadoriGosseGuerra}, see
  also~\cite{GuerraMarcelliniSchleper}.
\end{proof}

\section{Other Applications}
\label{sec:other}

The $p$-system~\eqref{eq:3} is of use in a variety of situations and
the procedure presented above may well be applied to them.

\subsection{Water Flowing in a Pipe}

\begin{figure}[h!]
  \centering 
\begin{picture}(0,0)%
\includegraphics{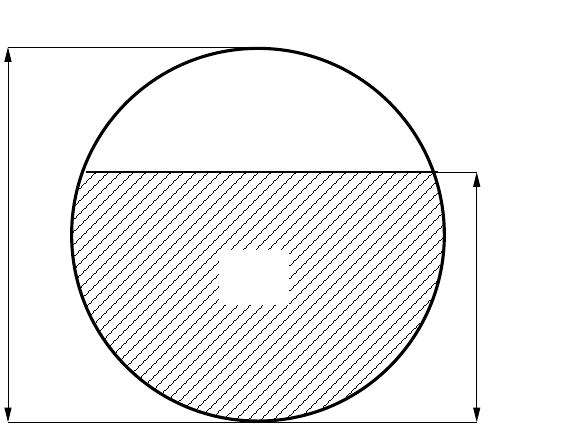}%
\end{picture}%
\setlength{\unitlength}{1973sp}%
\begingroup\makeatletter\ifx\SetFigFont\undefined%
\gdef\SetFigFont#1#2#3#4#5{%
  \reset@font\fontsize{#1}{#2pt}%
  \fontfamily{#3}\fontseries{#4}\fontshape{#5}%
  \selectfont}%
\fi\endgroup%
\begin{picture}(5477,4063)(2325,-5773)
\put(2476,-5311){\makebox(0,0)[lb]{\smash{{\SetFigFont{6}{7.2}{\familydefault}{\mddefault}{\updefault}{\color[rgb]{0,0,0}$d$}%
}}}}
\put(4576,-4411){\makebox(0,0)[lb]{\smash{{\SetFigFont{6}{7.2}{\familydefault}{\mddefault}{\updefault}{\color[rgb]{0,0,0}$A$}%
}}}}
\put(6601,-5386){\makebox(0,0)[lb]{\smash{{\SetFigFont{6}{7.2}{\familydefault}{\mddefault}{\updefault}{\color[rgb]{0,0,0}$h$}%
}}}}
\end{picture}%
\qquad
\begin{picture}(0,0)%
\includegraphics{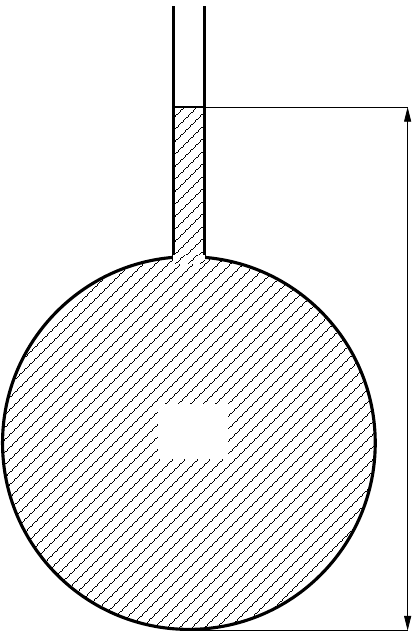}%
\end{picture}%
\setlength{\unitlength}{1973sp}%
\begingroup\makeatletter\ifx\SetFigFont\undefined%
\gdef\SetFigFont#1#2#3#4#5{%
  \reset@font\fontsize{#1}{#2pt}%
  \fontfamily{#3}\fontseries{#4}\fontshape{#5}%
  \selectfont}%
\fi\endgroup%
\begin{picture}(3954,6045)(2989,-5773)
\put(4576,-3886){\makebox(0,0)[lb]{\smash{{\SetFigFont{6}{7.2}{\familydefault}{\mddefault}{\updefault}{\color[rgb]{0,0,0}$A$}%
}}}}
\put(6526,-2086){\makebox(0,0)[lb]{\smash{{\SetFigFont{6}{7.2}{\familydefault}{\mddefault}{\updefault}{\color[rgb]{0,0,0}$h$}%
}}}}
\end{picture}%
  \caption{Notation used in~\eqref{eq:SaintVenant} and
    \eqref{eq:PressureLaw}. Left: the cross section of a standard pipe
    used in the modeling of free surface flows.  Right: a pipe with
    the fictitious Preissmann slot used to describe pressurized
    flows.}
  \label{fig:SingleTube}
\end{figure}
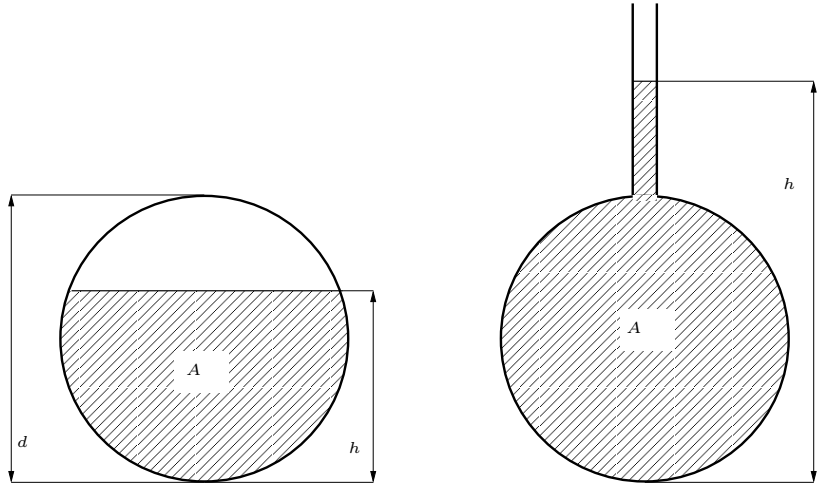

A different scenario that admits the same treatment presented in
Section~\ref{sec:m} is that of water flowing in a pipe. Neglecting
friction along the walls, in a horizontal pipe the Saint-Venant
equations~\cite{SaintVenant} read
\begin{equation}
  \label{eq:SaintVenant}
  \left\{
    \begin{array}{l}
      \partial_t a +\partial_x q = 0,
      \\[2mm]
      \partial_t q +\partial_x \left(\frac{q^2}{a}+p(a)\right) = 0 \,.
    \end{array}
  \right.
\end{equation}
Here, as usual, $t$ is time, $x$ the coordinate along the tube, $a =
a(t,x)$ is the area of the wet cross-section, $q = q(t,x)$ is the
water flow, so that $q = a \, v$, where $v = v(t,x)$ is the averaged
speed of water at time $t$ and position $x$. The hydrostatic term $p =
p(a)$ is defined as in~\cite[Section~3.2]{BorscheColomboGaravello},
namely
\begin{equation}
  \label{eq:PressureLaw}
  p(a) = g \, \int^a_0 \left(h(a) - h (\alpha)\right) \, \d\alpha
\end{equation}
where $h = h(a)$ is the height of water corresponding to $a$, see
Figure~\ref{fig:SingleTube}. Here $g$ is the acceleration due to
gravity.
In the case of water pipes, the function $h$ is often chosen
introducing the so-called Preissmann slot. It is an artificial
modification of the cross section of a tube, see
Figure~\ref{fig:SingleTube}, right, to merge free surface flow and
pressurized flow in a combined model. In the case of free surface flow
the physical geometry is used. In the case of pressurized flow a
narrow slot is added to the model, so that the \emph{height} of water
is extended beyond the tube diameter $d$. This widely used technique,
see, e.g.,~\cite{CungeHolly, Leon, Pagliara}, allows us to consider
both regimes in a single model.

With suitable choices of the term $f$, it is natural to consider the
following extension of~\eqref{eq:SaintVenant} to describe the dynamics
of water in a curved pipe:
\begin{equation}
  \label{eq:CurvedSaintVenant}
  \left\{
    \begin{array}{l}
      \partial_t a +\partial_x q = 0,
      \\[2mm]
      \partial_t q
      +
      \partial_x \left(\frac{q^2}{a}+p(a)\right)
      =
      - f (x) \, \kappa (x) \, q \,,
    \end{array}
  \right.
\end{equation}
with $p$ defined by the pressure law~\eqref{eq:PressureLaw}
satisfies~\textbf{(p)}. Referring to the case depicted in
Figure~\ref{fig:SingleTube} and
to~\cite[Section~3.2]{BorscheColomboGaravello}, calling $r$ the radius
of the pipe and $d$ the width of the Preissmann slot, we have
\begin{equation}
  \label{eq:h}
  h (a)
  =
      \begin{cases}
      \sqrt{\frac{2}{\pi} \, a},
      & a  \in  \left[0, \frac{\pi}{2}\, r^2\right] \,,
      \\[2mm]
      2r - \sqrt{2r^2 - \frac{2}{\pi} \, a},
      & a  \in
      (\frac{\pi}{2}\, r^2, \pi\, r^2 - \frac{1}{2\pi} \, d^2] \,,
      \\[2mm]
      \frac{a}{d} - \frac{1}{2\pi}\, d + 2r - \pi\, \frac{r^2}{d},
      & a  \in  (\pi r^2 - \frac{1}{2\pi} \, d^2, +\infty) \,.
    \end{cases}
\end{equation}
Straightforward computations show that the pressure
law~\eqref{eq:PressureLaw} with $h$ defined as in~\eqref{eq:h}
satisfies~\textbf{(p)}, so that the results in Section~\ref{sec:m} can
be applied also to~\eqref{eq:CurvedSaintVenant}.

\subsection{A Pipe with a Varying Section}

The dynamics of a fluid in a pipe with a slowly varying section $a = a
(x)$ is described by the well known equations
\begin{equation}
  \label{eq:7}
  \left\{
    \begin{array}{@{}l}
      \partial_t (a \, \rho) + \partial_x (a \, q) = 0,
      \\[2mm]
      \partial_t (a \, q)
      +
      \partial_x \! \left(a \left(\frac{q^2}{\rho} + p (\rho)\right)\right)
      = 0,
    \end{array}
  \right.
  \mbox{or }
  \left\{
    \begin{array}{@{}l@{}}
      \partial_t \rho + \partial_x q = -\frac{q}{a} \, \partial_x a,
      \\[2mm]
      \partial_t q
      +
      \partial_x \! \left(\frac{q^2}{\rho} + p (\rho)\right)
      =
      -\frac{q^2}{a\, \rho} \, \partial_x a,
    \end{array}
  \right.
\end{equation}
where $p = p (r)$ is the pressure law and, as in the previous section,
$\rho = \rho{(t,x)}$ is the fluid density and $q = q (t,x)$ is its
linear momentum density. The equivalence between the two
systems~\eqref{eq:7} is proved
in~\cite[Lemma~2.6]{ColomboMarcelliniJMAA}. This problem has been
widely considered in the literature, see, for
instance,~\cite{ChenHandbook, ColomboMarcelliniJMAA,
  DalMasoLeflochMurat, GoatinLefloch, KroenerLefloch, LiuTransonic}.

The system on the right in~\eqref{eq:7} clearly shows that a sudden
change in the pipe section, i.e., a discontinuity in the function $a$,
yields a Dirac delta function as source term in both
equations. Similarly to what was
\begin{figure}[!h]
  \centering \scalebox{0.2}{
\begin{picture}(0,0)%
\includegraphics{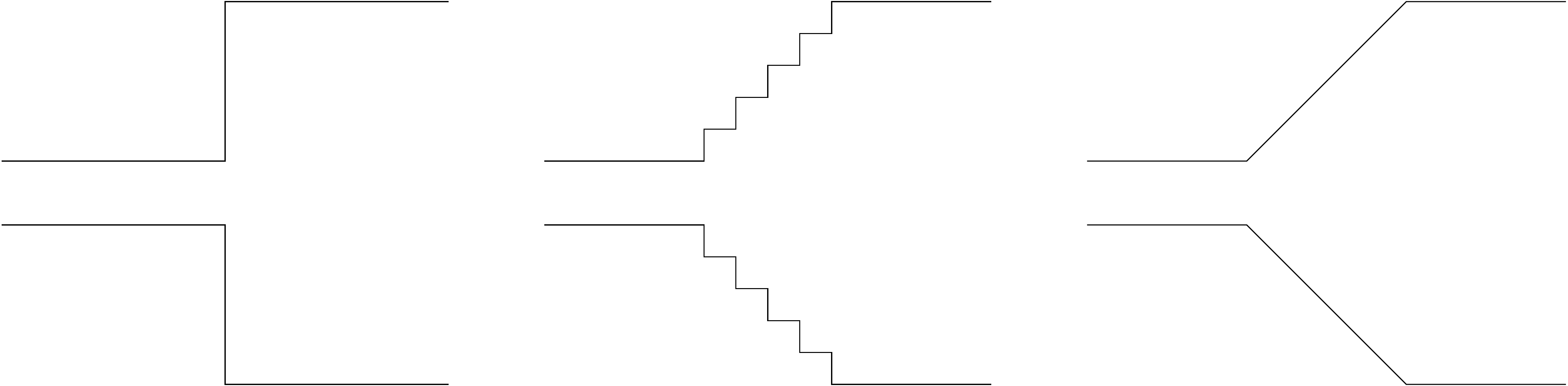}%
\end{picture}%
\setlength{\unitlength}{4144sp}%
\begingroup\makeatletter\ifx\SetFigFont\undefined%
\gdef\SetFigFont#1#2#3#4#5{%
  \reset@font\fontsize{#1}{#2pt}%
  \fontfamily{#3}\fontseries{#4}\fontshape{#5}%
  \selectfont}%
\fi\endgroup%
\begin{picture}(22094,5444)(-3621,-6833)
\end{picture}%
}
  \caption{Left: a single junction between two pipes. Middle: a
    sequence of junctions. Right: a pipe with a smoothly varying
    section.}
  \label{fig:section}
\end{figure}
done in Section~\ref{sec:m}, it is then natural to select a class of
solutions to~\eqref{eq:7} in the case of a single junction as in
Figure~\ref{fig:section}, left,
\begin{displaymath}
  a (x) =
    \begin{cases}
      a^-, & x  <  0 \,,
      \\
      a^+, & x  >  0 \,,
    \end{cases}
\end{displaymath}
pass to the case of a piecewise constant section $a = a (x)$ as in
Figure~\ref{fig:section}, center, and, in the limit, re-obtain
equations~\eqref{eq:7}. We refer to~\cite{ColomboMarcelliniJMAA} for
the details.

{\small

  \bibliographystyle{abbrv}


}

\end{document}